

   \magnification = \magstep 1
    \hsize = 6.3 true in
    \vsize = 9 true in
    \mathsurround 1pt
    \baselineskip = 20pt
			 \hoffset = -1.5mm

    \def\fine{$\bullet$}



    \font\smrm = cmr8
    \font\ssmrm = cmr6
     1


    2
   \font\tengothic=eufm10
   \font\sevengothic=eufm7
   \font\fivegothic=eufm5
   \font\tenmsx=msxm10
   \font\sevenmsx=msxm7
   \font\fivemsx=msxm5
   \font\tenmsy=msym10
   \font\sevenmsy=msym7
   \font\fivemsy=msym5

 \newfam\gothicfam
      \textfont\gothicfam=\tengothic
      \scriptfont\gothicfam=\sevengothic
      \scriptscriptfont\gothicfam=\fivegothic
 \newfam\msxfam
      \textfont\msxfam=\tenmsx
      \scriptfont\msxfam=\sevenmsx
      \scriptscriptfont\msxfam=\fivemsx
 \newfam\msyfam
      \textfont\msyfam=\tenmsy
      \scriptfont\msyfam=\sevenmsy
      \scriptscriptfont\msyfam=\fivemsy
   \def\goth#1{{\fam\gothicfam #1}}
   \def\bbb#1{{\fam\msyfam #1}}



   \def\Gr{Gr\"obner }
   \def\Definition{\noindent {\bf Definition. }}
   \def\Remark{\noindent {\bf Remark. }}

   \def\P{\bbb P}
			\def\N{\bbb N}
			\def\Z{\bbb Z}
   \def\Q{\bbb Q}
   
   \def\C{\bbb C}

		 \def\poly#1#2#3{\ifmmode#1\lbrack#2_1,#2_2,\dots,#2_{#3}\rbrack
			\else $#1\lbrack#2_1,#2_2,\dots,#2_{#3}\rbrack$\fi}


    \mathchardef\supsetneq="3A29
    \mathchardef\twoheadrightarrow="3910


   \def\proof{\noindent $\underline {Proof}$.\quad}


\def\authors#1 .{{\rm #1}}

\def\anno#1 .{{\rm #1}}

\def\title#1 .{{\sl#1}}

\def\journal#1 .{\frenchspacing #1}

\def\preprint#1 .{#1}

\def\book#1 .{{\it #1}}

\def\publisher#1 .{{\rm#1}}

\def\vol#1 .{{\bf #1}}

\def\pag#1 .{{\rm#1}}

\def\letterref#1 .{\hangindent=1.6cm\hangafter0\noindent
    \llap{\hbox to 1.6cm{\lbrack #1\rbrack \hskip 0cm plus 1cm}}}


\dimen3=\baselineskip \skip33= \dimen3 plus 6pt minus 6pt
\divide\skip33 by 4  \def\refskip{\vskip \skip33}


\def\MRefJour#1#2#3#4#5#6#7{\letterref#1 .\authors#2 .,
 \ \title#4 .,\ \journal#5 .\ \vol#6 . (\anno#3 .), \pag#7 .. \refskip}

\def\MRefBook#1#2#3#4#5{\letterref#1 .\authors#2 .,  \ \book#4 .,\
\publisher#5 .\ (\anno#3 .).\refskip}

\def\MRefInBook#1#2#3#4#5#6#7{\letterref#1 .\authors#2 .,  \ \title#4
.\ in \book#5 . \publisher#6 .\ (\anno#3 .), \pag#7 .. \refskip}

\def\MRefPrepr#1#2#3#4#5{\letterref#1 .\authors#2 .,  \ \title#4 .,\
(\anno#3 .), \ \preprint#5 .. \refskip}



\def\today{\ifcase\month\or January \or February \or March \or April
\or
   May \or June \or July \or August \or September \or October \or
November \or
   December \fi \space\number\day, \number\year}

\def\date{\parindent = 4truein Genova\quad \today \medskip}


\def\lef{\ssmrm \hfil G. Niesi, L. Robbiano \qquad}

\def\rig{\ssmrm Disproving Hibi's Conjecture \dots
\hfil}

\headline={\ifnum\pageno=1 \hfill \else \ifodd\pageno\rig
\else\lef \fi  \fi}

\vskip 3true cm
\noindent {\bf {Disproving Hibi's Conjecture with CoCoA}}

\noindent \hbox{\phantom{111111111111111111}} or

\noindent {\bf {Projective Curves with bad Hilbert Functions}}
\smallskip
\bigskip
\noindent Gianfranco Niesi, Lorenzo Robbiano  \footnote *{\smrm
The paper was partly supported by
Consiglio Nazionale delle Ricerche.}
\hfill \break
\noindent {\smrm Dipartimento di Matematica dell'Universit\'a di
Genova \ \  Via L.B. Alberti 4 \ \  16132 Genova \   ITALY }

\vskip 1 true cm

\medbreak
\beginsection Introduction.

In this paper we show how to combine different techniques from
Commutative Algebra and a systematic use of a Computer Algebra System
(in our case mainly CoCoA (see [G-N] and [A-G-N])) in order to
explicitly construct Cohen-Macaulay domains, which are standard
$k$-algebras and whose Hilbert function is \lq \lq bad".  In
particular we disprove a well-known conjecture by Hibi.

To be more precise, we recall that a ring $A$ is called a {\it
standard $k$-algebra}, or simply {\it standard}, if $k$ is a
field  and $A$ is a finitely generated $k$-algebra, which is
generated by its forms of degree $1$ (see [S$_1$]). To every
such a ring $A$ a numerical function $H_A$ is associated, namely
the function $H_A:\N \longrightarrow \N$, which is defined by
$H_A(r):=dim_kA_r$  for every $r\in \N$. Here it should be noted
that $A$ can be represented as the quotient of
$R:=k[X_0,\dots,X_n]$ modulo a homogeneous ideal $I$, hence
$dim_kA_r \leq dim_kR_r=  {{n+r }\choose {r }}$.  Such a
function is called the {\it Hilbert function} of $A$. It is
well-known (see [A-M]) that $H_A$ can be encoded in the power
series ${\cal P}_A(\lambda) :=\sum_rH_A(r)\lambda^r \in
\Z[[\lambda]]$, which is called the Hilbert-Poincar\'e series
(or simply the Poincar\'e series) of $A$. This series is
rational of type  ${\cal P}_A(\lambda) ={Q_A(\lambda) \over
(1-\lambda)^d}$, with $Q_A(1)\neq 0$; moreover
$Q_A(\lambda)=\sum h_i(A)\lambda^i \in \Z[\lambda]$ and if
$\delta$ is its degree, then the integral vector  ${\bf
h}(A):=(h_0(A),h_1(A),\dots,h_\delta(A))$ is called the {\it
$h$-vector of $A$}. It turns out that all the information of
$H_A$ can be encoded in $({\bf h}(A), d)$; in particular $d$ is
the dimension and $\sum_i h_i(A)=Q_A(1)$ is the multiplicity of
$A$. An efficient algorithm for the computation of ${\cal
P}_A(\lambda)$ is described in [B-C-R] and implemented in CoCoA.

\smallskip

\noindent In the paper [Hi] Hibi defines an $h$-vector
$(h_0,h_1,\dots,h_\delta)$ to be {\it flawless} if

i) $h_0 \leq h_1 \leq \cdots \leq h_{[\delta/2]}$ \quad and

ii) $h_i \leq h_{\delta-i}$ for every $i$ such that $0 \leq i
\leq [\delta/2]$

\noindent and he states the following

\noindent CONJECTURE: The $h$-vector ${\bf
h}(A):=(h_0(A),h_1(A),\dots,h_\delta(A))$ of a standard
Cohen-Macaulay domain is flawless.

The conjecture was supported by some empirical evidence and the
fact that the statement is true under some additional hypotheses
(see [Hi], Theorem 3.1).

The main goal of this paper is to construct {\it explicit
examples of standard Cohen-Macaulay domains,  whose $h$-vector
is not flawless}. To this end we use a technique introduced by
Galligo in [G], which yields sets of points in the projective
space with the Uniform Position Property. A good deal of freedom
in choosing equations allows us to construct projective
coordinate rings of sets of points, whose $h$-vector \lq \lq has
a flaw". Then we lift these sets to irreducible reduced rational
projective curves in $\P_{\C}^4$, which turn out to be
projectively Cohen-Macaulay. Their coordinate rings are the
desired counterexamples.

This paper is largely inspired by the work of Galligo (see [G])
and by the calculations of some Galois groups, which were shown
to the second author by G. Scheja, during a short visit to the
University of T\"ubingen in June 1991. To both we are largely
indebted.

\beginsection \S 1. Preliminaries.

In this section we recall all the definitions and results, which
we need later.

\medskip

\Definition Let $k$ be an infinite field, $A$ a standard
$k$-algebra, i.e. a graded $k$-algebra which is finitely
generated by its linear forms and let  ${\cal
P}_A(\lambda)={Q_A(\lambda) \over (1-\lambda)^d}$, with
$Q_A(1)\neq 0$, be its Poincar\'e series.  Let $\delta:=
deg(Q_A(\lambda))$ and
$Q_A(\lambda):=\sum_{i=0}^{\delta}h_i(A)\lambda^i$. Then ${\bf
h}(A):=(h_0(A),h_1(A),\dots,h_{\delta}(A))$ is called the
$h$-vector of $A$. Sometimes it is denoted by
$(h_0,h_1,\dots,h_{\delta})$, if there is no ambiguity about $A$.

 \medskip

The following facts are part of the folklore and are recalled
only for the sake of completeness. The non explained terminology
is part of the basic literature in Commutative Algebra and
Algebraic Geometry (see for instance [A-M] and [Hart]).

\medskip

\proclaim Lemma 1.1.  Let $k \subset F$ be fields, $A$ a
standard $k$-algebra, $A_F:=A\otimes_kF$. Then ${\bf h}(A_F)={\bf
h}(A)$.

\proof Indeed ${\cal P}_{ A_F}(\lambda)= {\cal P}_A(\lambda)$
\fine

 \medskip

\proclaim Proposition 1.2. Let k be an infinite field and A a
standard k-algebra of dimension d. Then there exist d linear
forms $L_1,\dots,L_d$, such that  $dim(A/(L_1,\dots,L_d))=0$. If
moreover A is Cohen-Macaulay (C-M), then $L_1,\dots,L_d$ is a
regular sequence in A.

\medskip

\proclaim Corollary 1.3. Let A be a Cohen-Macaulay standard
k-algebra over a field k, let $L_1,\dots,L_d$ be a maximal
regular sequence of linear forms in A and denote by $B:=
A/(L_1,\dots,L_d)$. Then ${\bf h}(A)={\bf h}(B)$.

\proof Indeed ${\cal P}_{B}(\lambda)= (1-\lambda)^d{\cal
P}_A(\lambda)$ \fine

\medskip

Let now $S:= k[X_0,X_1,\dots,X_n]$,  $\goth M$ a maximal
homogeneous relevant ideal in $S$ and $K:= K_0(S/{\goth M})$ its
associated field, i.e. the field of homogeneous fractions of
degree 0 of $S/{\goth M}$. The scheme $Proj(S/{\goth M})$ has a
unique point, whose associated local ring is $K_0(S/{\goth M})$.
If $X_0 \notin {\goth M}$, then we can dehomogenize ${\goth M}$
with respect to $X_0$ and we get a maximal ideal ${\goth m}$ in
$R:=k[X_1,\dots,X_n]$ (this can be done by putting $X_0=1$). It
turns out that  $K \cong R/{\goth m}$. Moreover every generic
linear change of coordinates yields the following shape of
${\goth m}$ $${\goth m}=(f(X_1),
X_2-g_2(X_1),X_3-g_3(X_1,X_2),\dots,
X_n-g_n(X_1,\dots,X_{n-1}))$$  It is clear that $K\cong
k[X]/(f(X))$, hence $deg(f(X) = dim_kK$.

\medskip \Definition In the above described situation we say that
$f(X)$ represents ${\goth M}$.

\medskip \Definition Let $\goth M$ be a maximal homogeneous
relevant ideal in $S$, $K$ its associated field and $d=dim_kK$ .
We say that $\goth M$ is $G$-symmetric if the Galois group
$Gal_k(K)$ is the full symmetric group $\Sigma_d$.

\medskip \proclaim Corollary 1.4. Let $\goth M$ be a maximal
homogeneous relevant ideal in $S$, $K$ its associated field and
$f(X)$ a polynomial of degree $d$ representing $\goth M$. Then
$\goth M$ is $G$-symmetric if and only if $Gal_k(f(X))=\Sigma_d$

\medskip Now we recall a well-known criterion   \medskip

\proclaim Theorem 1.5. Let $f(X) \in \Z[X]$, $d:=deg(f(X))$ and
assume that  \item {a)} $f(X)$ is irreducible  \item {b)} There
exist two prime numbers $p_1$, $p_2$ such that,  if we denote by
$f_i(X)$ the residue classes of $f(X)$ modulo $p_i$, $i=1,2$,
then  \itemitem {b$_1$)} $f_1(X)$ decomposes as the product of a
linear factor and an irreducible factor of degree $d-1$.
\itemitem {b$_2$)}  $f_2(X)$ decomposes as the product of an
irreducible factor of degree 2 and irreducible factors of odd
degrees. \hfill \break \noindent Then $Gal_{\Q}(f(X))= \Sigma_d$.

\proof See [W] and [S-S] \fine

\medskip

\Definition Let $F$ be an infinite field and let $E \subset
\P^n_F$ be a finite set of reduced points. Let $I$  be the
homogeneous ideal  of $F[X_0,\dots,X_n]$ which defines  $E$.
\quad We say that $E$ is  $G$-symmetric if there exists a
subfield $k \subseteq F$ and a maximal homogeneous relevant
ideal ${\goth M}$ in $S:=k[X_0,\dots,X_n]$ such that $I= {\goth
M}F[X_0,\dots,X_n]$ and  $Gal_k(K_0(S/{\goth M}))=\Sigma_d$

\proclaim Theorem 1.6. Let C be a projective irreducible reduced
curve in $\bbb P_{\C}^n$. Then the generic hyperplane section of
C is a G-symmetric set.

\proof See [G], Proposition 13. The proof given for a smooth
curve in $\bbb P^3$ works as well in general \fine

\proclaim Theorem 1.7. Let $E \subset \bbb P_F^n$ be a
G-symmetric set of points and $A:=F[X_0,\dots,X_n]/I$ its
coordinate ring. Let ${\bf h}(A) := (h_0,h_1,\dots,h_{\delta})$
be the h-vector of A. Then the following inequalities hold:
$$h_0+h_1+\cdots+h_i\leq h_{\delta-i}+\cdots+h_{\delta-1}+1$$
 for every $i=1,\dots, [\delta/2]$

\proof It is proved in [G] that if $E$ is $G$-symmetric, then it
has the Uniform Position Property (UPP), i.e. the Hilbert
function of its subsets depends only on their cardinality. Now
if $E$ is UPP, then every subset has the Cayley Bacharach (CB)
property (see [G-K-R] for definitions and properties), and the
conclusion follows again by [G-K-R] \fine

\beginsection \S 2. The construction of the counterexamples.

Now we are ready to use the above described machinery in order to
produce standard Cohen Macaulay domains with bad Hilbert
functions. All the following computations have been done using
CoCoA 1.7b on  a Macintosh.

The first step is to produce polynomials $f(X) \in \Z[X]$ of
degree $d$, such that $Gal_{\Q}(f(X)) = \Sigma_d$. This is not
difficult, since the "generic" one has this property; however we
want polynomials, which are not too dense, since we want to use
them to make further computations.

\medskip \proclaim Lemma 2.1. The polynomial $f(X):= X^{18}-X-1$
is such that $Gal_{\Q}(f(X)) = \Sigma_{18}$.

\proof We compute the factorization of $f(X)$ in $\Z[X]$  and
then the factorization of its classes modulo successive primes.
In this case we are particularly lucky, since we find that
\item{a)} $f(X)$ is irreducible  \itemitem{b$_1$)} The complete
factorization of $f(X)$ modulo 3 is
 $(X^2 - X - 1)(X^{13} - X^{12} + X^{11} + X^8 - X^7 - X^6 - X^5
+  X^4 + X^3 - X^2 + 1)(X^3 - X^2 + 1)$  \itemitem{b$_2$)} The
complete factorization of $f(X)$  modulo 5 is  $(X + 2)(X^{17} -
2X^{16} - X^{15} + 2X^{14} + X^{13} - 2X^{12} - X^{11} + 2X^{10}
+ X^9 - 2X^8 - X^7 + 2X^6 +  X^5 - 2X^4 - X^3 + 2X^2 + X + 2)$

\noindent The conclusion follows from Theorem 1.5 \fine

\medskip

\proclaim Corollary 2.2. Let I be the ideal of $\Q[X,Y,Z]$
generated by $(X^{18}-X-1, Y-g(X), Z-h(X,Y)$ with  $g(X) \in \bbb
Q[X]$ and $h(X,Y) \in \Q[X,Y]$. Let ${\goth M}:={}^hI$ i.e the
homogeneization of I with respect to a new indeterminate W. Then
${\goth M}$ is a G-symmetric maximal relevant ideal in $\bbb
Q[X,Y,Z,W]$.

\proof  Clearly $\Q[X,Y,Z]/I \cong \Q[X]/(X^{18}-X-1)$ and this
is a field, hence $I$ is a maximal ideal of $\Q[X,Y,Z]$.
Consequently ${\goth M}$ is a maximal relevant ideal in $\bbb
Q[X,Y,Z,W]$. The conclusion follows from Corollary 1.4 and Lemma
2.1 \fine

\medskip

Now the homogeneization of $I$ is computed via a \Gr basis
computation with respect to an ordering which is
degree-compatible and the leading term ideal of $I$ and of \
$^hI$ are generated by the same elements, hence they have the
same $h$-vector. The key point is now that we are {\it totally
free}\/ in the choice of $g(X)$ and $h(X,Y)$. We use again CoCoA
and again we are lucky, because we do not need many experiments.
Namely
 \medskip

\proclaim Example 2.3. Let I be the ideal of $\Q[X,Y,Z]$
generated by $(X^{18}-X-1, Y-X^3, Z-XY)$ and let ${\goth M}$  be
its homogeneization with respect to the new indeterminate W.
Then  $A:= \Q[X,Y,Z,W]/{\goth M}$ is a standard ${\Q}$-algebra,
which is a Cohen-Macaulay domain and whose $h$-vector is
(1,3,5,4,4,1). It satisfies the inequalities of Theorem 1.7, but
it is not flawless.

\proof Let $A:= \Q[X,Y,Z,W]/{\goth M}$. The computation  shows
that  ${\goth M}=(XY - ZW, \quad X^3 - YW^2, \quad  X^2Z -
Y^2W,  \quad Y^3 - XZ^2,  \quad Y^2Z^3 - XW^4 - W^5,  \quad
YZ^4 - X^2W^3 - XW^4,  \quad Z^5 -X^2W^3 - YW^4)$ and that ${\cal
P}_A(\lambda)={{(1 + 3\lambda + 5\lambda^2 + 4\lambda^3 +
4\lambda^4 + \lambda^5)} \over {(1-\lambda)}}$. Moreover $A$ is
a domain and it is Cohen-Macaulay, since it is 1-dimensional and
$W$ is a non zero divisor modulo ${\goth M}$ \fine

\medskip

This is already a counterexample to Hibi's conjecture!

However the fact that $A$ is a domain heavily relies on the
special ground field. Namely if we replace $\Q$ with $\C$ (it
suffices to replace it with the decomposition field of $f(X)$),
then $\C[X,Y,Z,W]/{\goth M} \C[X,Y,Z,W]$ is a reduced
$\C$-algebra, hence it is the coordinate ring of a $G$-symmetric
set of 18 points in $\bbb P_{\C}^3$. Its h-vector is still the
same (see Lemma 1.1), but  it is no more a domain.

So the final part is devoted to find a standard algebra whose
$h$-vector is not flawless and which is a \lq \lq geometric"
domain, i.e the fact that it is a domain is not affected by any
extension of the base field. The idea is to \lq \lq lift" our
previous example to an irreducible reduced curve in $\bbb
P^4_{\C}$. A small deformation of our equation $f(X)$ does the
trick. Namely

\medskip \proclaim Example 2.4.  Let ${\goth p}$ be the ideal of
\ $\C[X,Y,Z,T]$  generated by $(X^{18}-X-1-T, Y-X^3, Z-XY)$ and
let ${\goth P}$  be its homogeneization with respect to the new
indeterminate W. Then  $A:= \C[X,Y,Z,T,W]/{\goth P}$ is a
standard \ ${\C}$-algebra, which is a Cohen-Macaulay domain and
whose $h$-vector is (1,3,5,4,4,1), hence it is not flawless.

\proof It is clear that $\C[X,Y,Z,T]/{\goth p} \cong \C[X]$,
hence  ${\goth p}$ is a prime ideal, hence ${\goth P}$ is a
prime ideal. It defines a projective rational curve in
$\P_\C^4$. The actual computation yields the following minimal
system of generators for ${\goth P}$.

${\goth P}=(XY - ZW,\quad  X^3 - YW^2,\quad X^2Z - Y^2W, \quad
Y^3 - XZ^2, \quad Y^2Z^3 + XW^4 + TW^4 + W^5, \quad YZ^4 +
X^2W^3 + XTW^3 + XW^4, \quad Z^5 + X^2TW^2+ X^2W^3 + YW^4)$.

The computation of the Poincar\'e series yields ${\cal
P}_A(\lambda)= {{(1 + 3\lambda + 5\lambda^2 + 4\lambda^3 +
4\lambda^4 + \lambda^5)} \over {(1-\lambda)^2}}$. It remains to
prove that $A$ is a Cohen-Macaulay ring. For, it is enough to
show that $W,T$ is a regular sequence $mod \ {\goth P}$.  Indeed
$W$ is a non zero divisor $mod \ {\goth P}$, since it is the
homogeneizing indeterminate. If we compute the quotient modulo
$W$, we get $\C[X,Y,Z,T]/(XY,X^3,X^2Z,Y^3 -XZ^2,Y^2Z^3,YZ^4,Z^5
)$. Hence clearly $T$ does not divide zero \fine

\medskip

We conclude the paper with some remarks.

\medskip

\Remark Example 2.4 disproves Hibi's Conjecture. So it is
interesting to check that it does not fit with the special class
described by Hibi in [Hi] Theorem 3.1. There it is required that
the associated order ideal of monomials is pure. Without going
too much in to the details, we check that in our case the
associated order ideal of monomials is  $$\{1, X, Y, Z, X^2, XZ,
Y^2, YZ, Z^2, XZ^2, Y^2Z, YZ^2, Z^3, XZ^3, Y^2Z^2, YZ^3, Z^4,
XZ^4\}$$

\noindent Its maximal elements are $X^2, Y^2Z^2, YZ^3, XZ^4$
whose degrees are 2,4,4,5; therefore the order ideal of
monomials is  not pure.

\medskip

\Remark One can construct similar examples to 2.4. For instance
if we consider the ideals $(X^{22}-X-1-T, Y-X^3, Z-XY)$,
$(X^{26}-X-1-T, Y-X^3, Z-XY)$, $(X^{30}-X-1-T, Y-X^3, Z-XY)$ and
carry over the same construction as in 4.2, we see that:

the Galois group of $X^{22}-X-1$  is $\Sigma_{22}$. The primes
who do the trick as in Lemma 2.1 are 29 and 107. The
corresponding $h$-vector is (1,3,5,4,4,4,1).

The Galois group of $X^{26}-X-1$  is $\Sigma_{26}$. The primes
who do the trick as in Lemma 2.1 are 19 and 67. The corresponding
$h$-vector is (1,3,5,4,4,4,4,1).

The Galois group of $X^{30}-X-1$  is $\Sigma_{30}$. The primes
who do the trick as in Lemma 2.1 are 5 and 53. The corresponding
$h$-vector is (1,3,5,4,4,4,4,4,1).

All of them are not flawless.

\bigskip \noindent {\bf  References }  \smallskip

\MRefPrepr {A-G-N}{Armando, E., Giovini, A., Niesi, G.}{1991}
{CoCoA User's  Manual v.\ 1.5}{Dipartimento di
Ma\-te\-ma\-ti\-ca, Universit\`a di Genova}{}

\MRefBook {A-M}{Atiyah, M. F., Macdonald, I. G.}{1969}
{Introduction to Commutative Algebra} {Addison-Wesley}

\MRefJour {B-C-R}{Bigatti, A., Caboara, M., Robbiano,
L.}{1990}{On the computation of   Hilbert-Poincar\'e
series}{Applicable Algebra in Engineering, Communications and
Computing}{}{To appear}

\MRefJour {G}{Galligo, A.}{1991}{Exemples d'ensembles de Points
en Position Uniforme}{Proceedings of MEGA-90}{}{Birkhauser}

\MRefJour {G-K-R}{Geramita, A., Kreuzer, M., Robbiano, L.}{1991}
{Cayley-Bacharach schemes and their  canonical modules} {Trans.
Amer. Math. J.}{}{To appear}

\MRefJour{G-M-R}{Geramita, A.V., Maroscia, P., Roberts, L.G.}
{1983}{The Hilbert function of a reduced $k$-algebra}{J. London
Math. Soc.}{28}{443--452}

\MRefJour{G-M}{Geramita, A.V., Migliore, J.C.}{1989}{Hyperplane
sections of a smooth curve in ${\P^3}$}{Comm.
Algebra}{17}{3129--3164}

\MRefJour{G-N}{Giovini, A., Niesi, G.}{1990}{CoCoA: a
user-friendly system for commutative algebra}{Proceedings of
DISCO-90, Lecture Notes in Computer
Sciences}{429}{Springer-Verlag}

\MRefJour {Ha}{Harris, J.}{1980} {The genus of space
curves}{Math. Ann.}{249}{191--204}

\MRefInBook{H-E}{Harris, J. (with the collaboration of D.
Eisenbud)}{1982}{Curves in projective space}{S\'em.\ de
Math\'ematiques Sup\'erieures, Universit\'e de Montreal}{}{}

\MRefBook {Hart} {Hartshorne, R.}{1977}{Algebraic
Geometry}{Springer}

\MRefJour {Hi}{Hibi, T.}{1989}{Flawless O-sequences and Hilbert
functions of Cohen-Macaulay integral domains}{J. Pure Appl.
Algebra} {60}{245--251}

\MRefJour {Ra}{Rathmann, J.}{1987}{The uniform position principle
for curves in characteristic $p$}{Math. Ann.}{276}{565--579}

\MRefInBook {R$_1$}{Robbiano, L.}{1988}{Introduction to the
theory of Gr\"obner bases}{The Curves Seminar at Queen's, vol.\
V, Queen's Papers in Pure and Appl. Math.}{80}{Queen's
University, Kingston B1--B29}

\MRefJour {R$_2$} {Robbiano, L.}{1990}{On the theory of Hilbert
functions} {Queen's Papers in Pure and Appl. Math., Kingston
Canada}{85}{Vol VII}

\MRefBook {S-S}{Scheja, G., Storch, U.}{1988} {Lehrbuch der
Algebra} {B. G. Teubner Stuttgart}

\MRefJour {S$_1$}{Stanley, R}{1978} {Hilbert functions of graded
algebras} {Adv. in Math.}{28}{57--83}

\MRefPrepr{S$_2$}{Stanley, R}{1990}{On the Hilbert function of a
graded Cohen-Macaulay domain} {pre\-print, Massachusetts Inst.\
of Technology, Cambridge}

\MRefBook {W}{van der Waerden, B.L.}{1970}{Algebra, Vol I}{Ungar}

\bye